\documentstyle[11pt]{article}

\makeatletter

\@addtoreset{equation}{section}
\makeatother
\def\<{\langle}
\def\>{\rangle}

\newtheorem{lem}{Lemma}[section]
\newtheorem{theo}{Theorem}[section]
\newtheorem{rem}{Remark}[section]
\newtheorem{pro}{Proposition}[section]

\makeatletter
   
   \@addtoreset{equation}{section}
\makeatother

\begin{document}
\title{\bf Uniform Energy Decay for Wave Equations with Unbounded Damping Coefficients}
\author{Ryo IKEHATA\thanks{Corresponding author: ikehatar@hiroshima-u.ac.jp} \\ {\small Department of Mathematics, Graduate School of Education} \\{\small Hiroshima University} \\ {\small Higashi-Hiroshima 739-8524, Japan}\\ Hiroshi TAKEDA 
\thanks{h-takeda@fit.ac.jp} \\{\small Department of Intelligent Mechanical Engineering, Faculty of Engineering} \\{\small Fukuoka Institute of Technology} \\{\small Fukuoka 811-0295, Japan}}
\maketitle
\begin{abstract}
We consider the Cauchy problem for wave equations with unbounded damping coefficients in ${\bf R}^{n}$. For a general class of unbounded damping coefficients, we derive uniform total energy decay estimates together with a unique existence result of a weak solution. In this case we never impose strong assumptions such as compactness of the support of the initial data. This means that we never rely on the finite propagation speed property of the solutions, and we try to deal with an essential unbounded coefficient case. One of our methods comes from an idea developed in \cite{IM}. 
\end{abstract}
\section{Introduction}
\footnote[0]{Keywords and Phrases: Unbounded damping; Wave equation; Cauchy problem; Weighted initial data; Multiplier method; Fourier analysis; Total energy decay, Weak solutions.}
\footnote[0]{2010 Mathematics Subject Classification. Primary 35L05; Secondary 35B35, 35B40.}
We consider the mixed problem for wave equations with a localized damping in ${\bf R}^{n}$ ($n \geq 1$)
\begin{equation}
u_{tt}(t,x) -\Delta u(t,x) + a(x)u_{t}(t,x) = 0,\ \ \ (t,x)\in (0,\infty)\times {\bf R}^{n},\label{eqn}
\end{equation}
\begin{equation}
u(0,x)= u_{0}(x),\ \ u_{t}(0,x)= u_{1}(x),\ \ \ x\in{\bf R}^{n} ,\label{initial}
\end{equation}
where $(u_{0},u_{1})$ are initial data chosen as:
\[u_{0} \in H^{2}({\bf R}^{n}),\quad u_{1} \in H^{1}({\bf R}^{n}),\]
and
\[u_{t}=\frac{\partial u}{\partial t},\quad  u_{tt}=\frac{\partial^2 u}{\partial t^2}, \quad \Delta = \sum_{j=1}^{n}\frac{\partial^{2}}{\partial x_{j}^{2}}, \quad x = (x_1,\cdots,x_n).\]
Note that solutions and/or functions considered in this paper are all real valued except for several parts concerning the Fourier transform.

Concerning the decay or non-decay property of the total or local energy to problem (1.1)-(1.2) with $x$-dependent variable damping coefficients, many research manuscripts are already published by Alloui-Ibrahim-Khenissi \cite{AIK}, Bouclet-Royer \cite{BR}, Daoulatli \cite{D}, Ikehata \cite{Ik-1}, Ikehata-Todorova-Yordanov \cite{ITY}, Joly-Royer \cite{JR}, Kawashita \cite{K}, Khader \cite{Kha}, Matsumura \cite{Ma}, Mochizuki \cite{M}, Mochizuki-Nakazawa \cite{MN}, Nakao \cite{Na}, Nishiyama \cite{Nishiyama}, Nishihara \cite{Ni}, Radu-Todorova-Yordanov \cite{RTY}, Sobajima-Wakasugi \cite{SW-0}, Todorova-Yordanov \cite{TY}, Uesaka \cite{U} and Wakasugi \cite{W}, Zhang \cite{Z} and the references therein. However, we should emphasize that those cases are quite restricted to the bounded damping coefficient case, i.e., $a \in L^{\infty}({\bf R}^{n})$. This condition seems to be essential to get the unique existence of mild or weak solutions to problem (1.1)-(1.2). So, when we do not assume the boundedness of the coefficient $a(x)$, a natural question arises whether one can construct a unique weak or mild solution $u(t,x)$ to problem (1.1)-(1.2) together with some decay property of the total energy or not. Quite recently Sobajima-Wakasugi \cite{SW} have announced an interesting result from the viewpoint of the diffusion phenomenon of the solution to the equation (1.1) together with a unique global existence result. The mixed problem to the equation (1.1) in \cite{SW} is considered in the exterior domain of a bounded obstacle, and the Dirichlet null boundary condition is treated. They treated typically $a(x) = a_{0}\vert x\vert^{\alpha}$ with $a_{0} > 0$ and $\alpha > 0$. But, their results require a stronger assumption such as the compactness of the support of the initial data. This implies that they have to rely on the finite speed of propagation property (FSPP for short) of the solution, so that in an essential meaning, their framework seems to be still restricted to the bounded damping coefficient case for each $t \geq 0$. In the treatment of the unbounded coefficient $a(x)$, it seems important and interesting not to assume such FSPP. Additionally, they essentially used rather stronger regularity condition on the initial data such that $[u_{0},u_{1}] \in (H^{2}\cap H_{0}^{1})\times H_{0}^{1}$. In connection with this topic, D'Abbicco \cite{Da} (and for a more general class, see D'Abbicco-Ebert \cite{DE} and Reissig \cite{R}) and Wirth \cite{Wi} have ever studied $t$-dependent unbounded damping coefficient case:
\[u_{tt}(t,x) -\Delta u(t,x) + b(1+t)^{\alpha}u_{t}(t,x) = f(u),\]       
where $\alpha \in [0,1]$ and $b > 0$. Therefore, in the $x$-dependent unbounded coefficient case, a unique existence of the solution itself together with some decay property of the total energy are completely open in the framework of non-compactly supported initial data class. When Sobajima-Wakasugi \cite{SW} constructs a unique solution, they have used directly the well-known result due to Ikawa \cite{I}. That means their result is still under the already known framework. In our result to be announced, we have to discuss how we should construct a weak solution itself because no any well-known theories can be applied directly. This is a main difficulty in our result.

In this connection, originally Komatsu \cite{Ko} first proposed this open question in his Master thesis in January 2016 such that for the unbounded damping coefficient case $a \notin L^{\infty}({\bf R}^{n})$, can one construct a global in time solution?\, Unfortunately, Komatsu \cite{Ko} could not solve his problem before finishing Master course. This paper gives an answer to his problem.      

Now, let us start with introducing our new result. Before stating our result, we shall give the following only one assumption ({\bf A}) on the damping coefficient $a(x)$, and the definition of the solution to be constructed.\\

\noindent
({\bf A})\,$a \in C({\bf R}^{n})$, and there exists a constant $V_{0} > 0$ such that $0 < V_{0} \leq a(x)$ ($\forall x \in {\bf R}^{n}$).\\

\noindent
{\bf Definition.} {\rm A function $u:\,[0,\infty)\times{\bf R}^{n} \to {\bf R}$ is called as the weak solution if it satisfies
\[\int_{0}^{\infty}\int_{{\bf R}^{n}}u(t,x)(\phi_{tt}(t,x)-\Delta\phi(t,x) - a(x)\phi_{t}(t,x))dx dt\] 
\[= \int_{{\bf R}^{n}}u_{1}(x)\phi(0,x)dx - \int_{{\bf R}^{n}}u_{0}(x)\phi_{t}(0,x)dx+ \int_{{\bf R}^{n}}a(x)u_{0}(x)\phi(0,x)dx\]
for any $\phi \in C_{0}^{\infty}([0,\infty)\times{\bf R}^{n})$.}\\

Our new result reads as follows.  

\begin{theo} Let $n \geq 3$ and assume {\rm ({\bf A})}. If the initial data $[u_{0},u_{1}] \in (H^{2}({\bf R}^{n})\cap L^{1}({\bf R}^{n}))\times (H^{1}({\bf R}^{n})\cap L^{1}({\bf R}^{n}))$ further satisfies $a(\cdot)u_{0} \in L^{1}({\bf R}^{n})\cap L^{2}({\bf R}^{n})$, then there exists a unique weak solution $u \in L^{\infty}(0,\infty;H^{1}({\bf R}^{n}))\cap W^{1,\infty}(0,\infty;L^{2}({\bf R}^{n}))$ to problem {\rm (1.1)}-{\rm (1.2)} satisfying
\[\Vert u(t,\cdot)\Vert^{2} \leq CI_{00}^{2}(1+t)^{-1},\quad E_{u}(t) \leq CI_{00}^{2}(1+t)^{-2},\]
with some constant $C > 0$, where
\[I_{00} := \left(\Vert u_{0}\Vert^{2} + \Vert\nabla u_{0}\Vert^{2} + \Vert a(\cdot)u_{0}\Vert_{1}^{2} + \Vert a(\cdot)u_{0}\Vert^{2} + \Vert u_{1}\Vert^{2} +\Vert u_{1}\Vert_{1}^{2}\right)^{1/2}.\]
\end{theo}

\begin{rem}{\rm If one considers the mixed problem (1.1)-(1.2) with the Dirichlet null boundary condition on the smooth exterior domain, one can treat the two dimensional case. In that case, we need to assume the logarithmic type weight condition such that $\Vert\log(B\vert x\vert)(a(\cdot)u_{0} + u_{1})\Vert < +\infty$ on the initial data with some constant $B > 0$. This has a close relation to the place where one obtains Lemma 2.1 below in the two dimensional exterior domain case. For more detail, one can refer the reader to \cite{IM}.}
\end{rem}
\begin{rem}{\rm The assumption $[u_{0},u_{1}] \in H^{2}({\bf R}^{n})\times H^{1}({\bf R}^{n})$ is not essential. It will be used only to justify the integration by parts in the course of the proof. That condition seems not to be so rare (cf. \cite{MN}). One may be able to generalize to the more weak case such that $[u_{0},u_{1}] \in H^{1}({\bf R}^{n})\times L^{2}({\bf R}^{n})$, however, for simplicity we do not deal with such case.}
\end{rem}
\noindent
{\bf Example.} {\rm As the typical unbounded example for $a(x)$, we can choose $a(x) := (1+\vert x\vert^{2})^{\frac{\alpha}{2}}$ with $\alpha \in [0,\infty)$, $a(x) := e^{\vert x\vert}$, and so on.}
\par
\vspace{0.2cm}

This paper is organized as follows. In section 2 we shall prove Theorem 1.1 by relying on a multiplier method which was introduced in \cite{IM}. In section 3 we give several remarks and open problems. Section 4 is devoted to the appendix to check the known result.\\

{\bf Notation.} {\small Throughout this paper, $\| \cdot\|_q$ stands for the usual $L^q({\bf R}^{n})$-norm. For simplicity of notation, in particular, we use $\| \cdot\|$ instead of $\| \cdot\|_2$. Furthermore, we denote $\Vert\cdot\Vert_{H^{1}}$ as the usual $H^{1}$-norm. The $L^{2}$-inner product is denoted by $(f,g) := \displaystyle{\int_{{\bf R}^{n}}}f(x)g(x)dx$ for $f,g \in L^{2}({\bf R}^{n})$. The total energy $E_{u}(t)$ corresponding to the solution $u(t,x)$ of (1.1) is defined by
\[E_{u}(t):=\frac{1}{2}(\| u_t(t,\cdot)\|^2+\|\nabla u(t,\cdot)\|^2),\]
where 
\[\vert\nabla f(x)\vert^{2} := \sum_{j=1}^{n}\vert\frac{\partial f(x)}{\partial x_{j}}\vert^{2}.\]
The weighted $L^{1}$-space and its norm $\Vert \cdot\Vert_{1,m}$ can be defined as
\[f \in L^{1,m}({\bf R}^{n}) \Leftrightarrow f \in L^{1}({\bf R}^{n}),\quad \Vert f\Vert_{1,m} := \int_{{\bf R}^{n}}(1+\vert x\vert)^{m}\vert f(x)\vert dx < +\infty.\]
The subspace $X_{m}^{n}$ of $L^{1,m}$ is defined by
\[X_{m}^{n} := \{f \in L^{1,m}({\bf R}^{n}):\,\int_{{\bf R}^{n}}f(x)dx = 0\}.\]
And also, for the Hilbert space $X$ we define a class of vector valued continuous functions $C_{0}([0,\infty);X)$ as follows: $f \in C_{0}([0,\infty);X)$ if and only if $f \in C([0,\infty);X)$ and the closure of the set $\{t \in [0,\infty);\,\Vert f(t)\Vert_{X} \ne 0\}$ is compact in $[0,\infty)$.

On the other hand, we denote the Fourier transform of $f(x)$ by $\hat{f}(\xi) := (\displaystyle{\frac{1}{2\pi}})^{\frac{n}{2}}\displaystyle{\int_{{\bf R}^{n}}}e^{-ix\cdot\xi}f(x)dx$ as usual with $i := \sqrt{-1}$, and we define the usual convolution by 
\[(f*g)(x) := \int_{{\bf R}^{n}}f(x-y)g(y)dy.\]
}\\ 

\section{Proof of Theorem 1.1.}

In the course of the proof, the next inequality concerning the Fourier image of the Riesz potential plays an alternative role for the Hardy inequality. This comes from \cite[Proposition 2.1]{Ike-0}. It should be emphasized that this inequality holds even in the low dimensional case if we control the weight parameter $\gamma$. 
\begin{pro}\,Let $n \geq 1$ and $\gamma \in [0,1]$.\\ 
{\rm (1)}\, If $f \in L^{2}({\bf R}^{n})\cap L^{1,\gamma}({\bf R}^{n})$, and $\theta \in [0,\displaystyle{\frac{n}{2}})$, then there exists a constant $C = C_{n,\theta,\gamma} > 0$ such that\\
\[\displaystyle{\int_{{\bf R}^{n}}}\displaystyle{\frac{\vert \hat{f}(\xi)\vert^{2}}{\vert\xi\vert^{2\theta}}}d\xi \leq C\left(\Vert f\Vert_{1,\gamma}^{2} + \vert\displaystyle{\int_{{\bf R}^{n}}f(x)dx}\vert^{2}+ \Vert f\Vert^{2}\right) .\]
{\rm (2)}\, If $f \in L^{2}({\bf R}^{n})\cap X_{\gamma}^{n}$, and $\theta \in [0,\gamma + \displaystyle{\frac{n}{2}})$, then it is true that 
\[\displaystyle{\int_{{\bf R}^{n}}}\displaystyle{\frac{\vert \hat{f}(\xi)\vert^{2}}{\vert\xi\vert^{2\theta}}}d\xi \leq C(\Vert f\Vert_{1,\gamma}^{2} + \Vert f\Vert^{2}) \]
with some constant $C = C_{n,\theta,\gamma} > 0$.
\end{pro}
 
To construct a global weak solution we first define a sequence of the weak solutions $\{u^{(m)}(t,x)\}$ ($m \in {\bf N}$) to the approximated problem below:
\begin{equation}
u_{tt}^{(m)}(t,x) -\Delta u^{(m)}(t,x) + a_{m}(x)u_{t}^{(m)}(t,x) = 0,\ \ \ (t,x)\in (0,\infty)\times {\bf R}^{n},\label{eqn2}
\end{equation}
\begin{equation}
u^{(m)}(0,x)= u_{0}(x),\ \ u_{t}^{(m)}(0,x)= u_{1}(x),\ \ \ x\in{\bf R}^{n},\label{initial2}
\end{equation}
where $a_{m} \in C({\bf R}^{n})$ can be chosen to satisfy
\begin{eqnarray}a_{m}(x) = \left\{
  \begin{array}{ll}
   \displaystyle{a(x)}&
       \qquad (\vert x\vert \leq m) \\[0.2cm]
   \displaystyle{V_{0}}& \qquad (\vert x\vert > m+1),
   \end{array} \right. \end{eqnarray}
and 
\begin{equation}
V_{0} \leq a_{m}(x) \leq a(x), \quad a_{m}(x) \to a(x) \quad \textstyle{as}\quad m \to \infty \quad\textstyle{(pointwise)}\quad x \in {\bf R}^{n}
\end{equation}
for each $x \in {\bf R}^{n}$.\\

Now let us consider the problem (2.1)-(2.2) with initial data $[u_{0}, u_{1}] \in H^{2}({\bf R}^{n})\times H^{1}({\bf R}^{n})$. Then, for each $m \in {\bf N}$ since $a_{m} \in C({\bf R}^{n})\cap L^{\infty}({\bf R}^{n})$ it is well known that the Cauchy problem (2.1)-(2.2) has a unique strong solution $u^{(m)} \in C([0,\infty);H^{2}({\bf R}^{n}))\cap C^{1}([0,\infty);H^{1}({\bf R}^{n}))\cap C^{2}([0,\infty);L^{2}({\bf R}^{n}))$ satisfying the energy identity:
\begin{equation}
E_{u^{(m)}}(t) + \int_{0}^{t}\int_{{\bf R}^{n}}a_{m}(x)\vert u_{s}^{(m)}(s,x)\vert^{2}dxds = E_{}(0),
\end{equation}
where
\[E(0) := \frac{1}{2}(\Vert u_{1}\Vert^{2} + \Vert\nabla u_{0}\Vert^{2}).\]

To begin with we prove the following crucial estimate. The lemma below is a combination of the method introduced in \cite{IM} ($=$ the modified Morawetz method \cite{Mora}) and Proposition 2.1.
\begin{lem}\,Let $n \geq 3$. Under the same assumptions as in Theorem {\rm 1.1}, the {\rm (}unique{\rm )} solution $u^{(m)}(t,x)$ to problem {\rm (2.1)-(2.2)} satisfies
\[\Vert u^{(m)}(t,\cdot)\Vert^{2} + \int_{0}^{t}\int_{{\bf R}^{n}}a_{m}(x)\vert u^{(m)}(s,x)\vert^{2}dxds\]
\[\leq C\left(\Vert a(\cdot)u_{0}\Vert_{1}^{2} + \Vert a(\cdot)u_{0}\Vert^{2} + \Vert u_{1}\Vert_{1}^{2} + \Vert u_{1}\Vert^{2}\right) =:CI_{0}^{2}\]
with a constant $C > 0$, where $C$ is independent of $m$.
\end{lem}
{\it Proof.}\, The original idea comes from \cite{IM}. For the solution $u^{(m)}(t,x)$ to problem (2.1)-(2.2), one introduces an auxiliary function
$$W(t,x) := \int^{t}_{0}u^{(m)}(s,x)ds.$$
Then $W(t,x)$ satisfies
\begin{equation}
W_{tt} - \Delta W + a_{m}(x)W_{t} = a_{m}(x)u_{0} + u_{1},\ \ \ \ (t,x) \in (0,\infty) \times {\bf R}^{n},
\end{equation}
\begin{equation}
W(0,x) = 0,\quad W_{t}(0,x) = u_{0}(x),\,\,\,\, x \in {\bf R}^{n}. 
\end{equation}
Multiplying $(2.6)$ by $W_{t}$ and integrating over $[0,t]\times {\bf R}^{n}$ we get
$$\frac{1}{2}(\|W_{t}(t,\cdot)\|^{2} + \|\nabla W(t,\cdot)\|^{2}) + \int^{t}_{0}\|\sqrt{a_{m}(\cdot)}W_{s}(s,\cdot)\|^{2}ds$$
\begin{equation}
= \frac{1}{2}\|u_{0}\|^{2} + \int^{t}_{0}(a_{m}(\cdot)u_{0} + u_{1}, W_{s}(s,\cdot))ds.
\end{equation}

Next one uses (1) of Proposition 2.1 with $\theta = 1$ and $\gamma = 0$, the Plancherel theorem and the Cauchy-Schwarz inequality to obtain a series of inequalities below:\\

$$\left\vert\int^{t}_{0}(a_{m}(\cdot)u_{0} + u_{1},W_{s}(s,\cdot))ds\right\vert = \left\vert\int^{t}_{0}\frac{d}{ds}(a_{m}(\cdot)u_{0} + u_{1},W(s,\cdot))ds\right\vert$$
$$= \left\vert\int_{{\bf R}^{n}}(a_{m}(x)u_{0}(x) + u_{1}(x))W(t,x)dx\right\vert$$
$$= \left\vert\int_{{\bf R}_{\xi}^{n}}(\widehat{(a_{m}u_{0})}(\xi) + \hat{u}_{1}(\xi))\overline{\hat{W}(t,\xi)}d\xi\right\vert$$
$$\leq \int_{{\bf R}_{\xi}^{n}}\frac{\vert\widehat{(a_{m}u_{0})}(\xi) + \hat{u}_{1}(\xi)\vert}{\vert\xi\vert}(\vert\xi\vert\vert\hat{W}(t,\xi)\vert)d\xi$$
$$\leq \left(\int_{{\bf R}_{\xi}^{n}}\frac{\vert\widehat{(a_{m}u_{0})}(\xi) + \hat{u}_{1}(\xi)\vert^{2}}{\vert\xi\vert^{2}}d\xi\right)^{1/2}\left(\int_{{\bf R}_{\xi}^{n}}\vert\xi\vert^{2}\vert\hat{W}(t,\xi)\vert^{2}d\xi\right)^{1/2}$$
\[\leq \int_{{\bf R}_{\xi}^{n}}\frac{\vert\widehat{(a_{m}u_{0})}(\xi) + \hat{u}_{1}(\xi)\vert^{2}}{\vert\xi\vert^{2}}d\xi + \frac{1}{4}\int_{{\bf R}_{\xi}^{n}}\vert\xi\vert^{2}\vert\hat{W}(t,\xi)\vert^{2}d\xi\]
\[\leq C\Vert a_{m}u_{0}+u_{1}\Vert_{1}^{2} + \Vert a_{m}u_{0}+u_{1}\Vert^{2} + C\left\vert\int_{{\bf R}^{n}}(a_{m}(x)u_{0}(x)+u_{1}(x))dx\right\vert^{2} + \frac{1}{4}\Vert\nabla W(t,\cdot)\Vert^{2}\]
\begin{equation}
\leq C\left(\Vert a_{m}(\cdot)u_{0}\Vert_{1}^{2} +\Vert u_{1}\Vert_{1}^{2} + \Vert a_{m}(\cdot)u_{0}\Vert^{2} +\Vert u_{1}\Vert^{2} + \vert\int_{{\bf R}^{n}}(a_{m}(x)u_{0}(x)+u_{1}(x))dx\vert^{2}\right) + \frac{1}{4}\Vert\nabla W(t,\cdot)\Vert^{2}
\end{equation}
with some constant $C > 0$. 
Combining $(2.8)$ with $(2.9)$ we can derive
$$\frac{1}{2}\|W_{t}(t,\cdot)\|^{2} + \frac{1}{4}\|\nabla W(t,\cdot)\|^{2} + \int^{t}_{0}\int_{{\bf R}^{n}}a_{m}(x)\left\vert W_{s}(s,x)\right\vert^{2}dxds$$
$$\leq C\left(\Vert a_{m}(\cdot)u_{0}\Vert_{1}^{2} +\Vert u_{1}\Vert_{1}^{2} + \Vert a_{m}(\cdot)u_{0}\Vert^{2} +\Vert u_{1}\Vert^{2}+ \vert \int_{{\bf R}^{n}}(a_{m}(x)u_{0}(x)+u_{1}(x))dx\vert^{2}\right)$$
with some constant $C > 0$. Since $a_{m}(x) \leq a(x)$, in the case when $n \geq 3$ one has
$$\frac{1}{2}\|W_{t}(t,\cdot)\|^{2} + \frac{1}{4}\|\nabla W(t,\cdot)\|^{2} + \int^{t}_{0}\int_{{\bf R}^{n}}a_{m}(x)\left\vert W_{s}(s,x)\right\vert^{2}dxds$$
$$\leq C\left(\Vert a(\cdot)u_{0}\Vert_{1}^{2} +\Vert u_{1}\Vert_{1}^{2} + \Vert a(\cdot)u_{0}\Vert^{2} +\Vert u_{1}\Vert^{2} + \vert\int_{{\bf R}^{n}}(a(x)\vert u_{0}(x)\vert + \vert u_{1}(x)\vert)dx\vert^{2}\right).$$
We easily see that the constant $C$ in the above estimates is independent of $m$.
Thus, one has the desired estimate because of the fact $W_{t}=u^{(m)}$.
$\hfill\Box$

\begin{lem}Under the same assumptions as in Theorem {\rm 1.1},  the {\rm (}unique{\rm )} solution $u^{(m)}(t,x)$ to problem {\rm (2.1)-(2.2)} satisfies
\[(1+t)E_{u^{(m)}}(t) \leq E(0)(1+\frac{1}{V_{0}} + \frac{1}{2\varepsilon}) + \frac{1}{2}(u_{1},u_{0}) =: I_{1}^{2},\quad (t \geq 0),\]
\[\int_{0}^{t}E_{u^{(m)}}(s)ds \leq I_{1}^{2},\quad (t \geq 0),\]
with some small constant $\varepsilon > 0$, where $C$ is independent of $m$.
\end{lem}
{\it Proof.}\,For simplicity of notation, we use $w(t,x)$ in place of $u^{(m)}(t,x)$, i.e., $w(t,x)$ satisfies the equation below:
\begin{equation}
w_{tt}(t,x) -\Delta w(t,x) + a_{m}(x)w_{t}(t,x) = 0,\ \ \ (t,x)\in (0,\infty)\times {\bf R}^{n},
\end{equation}
\begin{equation}
w(0,x)= u_{0}(x),\ \ w_{t}(0,x)= u_{1}(x),\ \ \ x\in{\bf R}^{n}.
\end{equation}
Note that $E_{w}(t) = E_{u^{(m)}}(t)$ satisfies (2.5). Then, since
\[\frac{d}{dt}\{(1+t)E_{w}(t)\} \leq E_{w}(t),\]
it follows from (2.5) that
\[(1+t)E_{w}(t) \leq E(0) + \frac{1}{2}\int_{0}^{t}\Vert w_{s}(s,\cdot)\Vert^{2}ds + \frac{1}{2}\int_{0}^{t}\Vert\nabla w(s,\cdot)\Vert^{2}ds\]
\[\leq E(0) + \frac{1}{2V_{0}}\int_{0}^{t}\int_{{\bf R}^{n}}a_{m}(x)\vert w_{s}(s,x)\vert^{2}dxds + \frac{1}{2}\int_{0}^{t}\Vert\nabla w(s,\cdot)\Vert^{2}ds\]
\begin{equation}
\leq E(0) + \frac{1}{2V_{0}}E(0) + \frac{1}{2}\int_{0}^{t}\Vert\nabla w(s,\cdot)\Vert^{2}ds.
\end{equation}
On the other hand, by multiplying both sides of (2.10) by $w(t,x)$ it follows that
\begin{equation}
\frac{d}{dt}(w_{t}(t,\cdot),w(t,\cdot)) + \Vert\nabla w(t,\cdot)\Vert^{2} + \frac{1}{2}\frac{d}{dt}\int_{{\bf R}^{n}}a_{m}(x)\vert w(t,x)\vert^{2}dx = \Vert w_{t}(t,\cdot)\Vert^{2},
\end{equation}
so that by integrating both sides over $[0,t]$ one has
\[\int_{0}^{t}\Vert\nabla w(s,\cdot)\Vert^{2}ds + \frac{1}{2}\int_{{\bf R}^{n}}a_{m}(x)\vert w(t,x)\vert^{2}dx\]
\[= \int_{0}^{t}\Vert w_{s}(s,\cdot)\Vert^{2}ds - (w_{t}(t,\cdot),w(t,\cdot)) + (u_{1},u_{0})\]
\[\leq \frac{1}{V_{0}}\int_{0}^{t}\int_{{\bf R}^{n}}a_{m}(x)\vert w_{s}(s,x)\vert^{2}dxds +\frac{1}{2\varepsilon}\Vert w_{t}(t,\cdot)\Vert^{2} + \frac{\varepsilon}{2}\Vert w(t,\cdot)\Vert^{2} + (u_{1},u_{0})\]
\[\leq \frac{1}{V_{0}}\int_{0}^{t}\int_{{\bf R}^{n}}a_{m}(x)\vert w_{s}(s,x)\vert^{2}dxds +\frac{1}{\varepsilon}E_{w}(t) + \frac{\varepsilon}{2V_{0}}\int_{{\bf R}^{n}}a_{m}(x)\vert w(t,x)\vert^{2}dx + (u_{1},u_{0})\]
\[\leq \frac{1}{V_{0}}E(0) + \frac{1}{\varepsilon}E(0) + \frac{\varepsilon}{2V_{0}}\int_{{\bf R}^{n}}a_{m}(x)\vert w(t,x)\vert^{2}dx + (u_{1},u_{0}),\]
where we have just used (2.5) and the Cauchy-Schwarz inequality with some positive parameter $\varepsilon > 0$.  This implies
\[\int_{0}^{t}\Vert\nabla w(s,\cdot)\Vert^{2}ds + \frac{1}{2}(1-\frac{\varepsilon}{V_{0}})\int_{{\bf R}^{n}}a_{m}(x)\vert w(t,x)\vert^{2}dx\]
\begin{equation}
\leq E(0)(\frac{1}{V_{0}} + \frac{1}{\varepsilon}) + (u_{1},u_{0}).
\end{equation}
By choosing $\varepsilon > 0$ sufficiently small, from (2.13) and (2.14) one can get the desired two estimates. In this final check, because of (2.6) we have to make the following estimate once more:
\[\int_{0}^{t}\Vert w_{t}(s,\cdot)\Vert^{2}ds \leq \frac{1}{V_{0}}\int_{{\bf R}^{n}}a_{m}(x)\vert w_{s}(s,x)\vert^{2}dxds \leq \frac{1}{V_{0}}E(0).\]
$\hfill\Box$

\begin{lem}Under the same assumptions as in Theorem {\rm 1.1},  the {\rm (}unique{\rm )} solution $u^{(m)}(t,x)$ to problem {\rm (2.1)-(2.2)} satisfies
\[(1+t)^{2}E_{u^{(m)}}(t) \leq E(0) + \frac{1}{V_{0}}(E(0) + I_{1}^{2}) + (u_{1},u_{0}) + \frac{1}{2}\int_{{\bf R}^{n}}a(x)\vert u_{0}(x)\vert^{2}dx + \frac{C}{2}I_{0}^{2} +\frac{I_{1}^{2}}{\varepsilon}=: I_{2}^{2},\]
\[(1+t)\Vert u^{(m)}(t,\cdot)\Vert^{2} \leq 2(V_{0}-\epsilon)^{-1}\{(u_{1},u_{0}) + \frac{1}{2}\int_{{\bf R}^{n}}a(x)\vert u_{0}(x)\vert^{2}dx + \frac{C}{2}I_{0}^{2} + \frac{1}{V_{0}}(E(0) + I_{1}^{2})+\frac{I_{1}^{2}}{\varepsilon} \} =: I_{3}^{2}\]
with some constant $C > 0$, where $C$ is independent of $m$.
\end{lem}
{\it Proof.}\,We use the same notation as in the proof of Lemma 2.2. Then, by multiplying both sides of (2.11) by $(1+t)w$, and integrating it over $[0,t] \times {\bf{R}}^{n}$ one can arrive at
the important identity:
\[\frac{1}{2}\Vert u_{0}\Vert^{2} + \int_{0}^{t}(1+s)\Vert\nabla w(s,\cdot)\Vert^{2}ds + \frac{(1+t)}{2}\int_{{\bf R}^{n}}a_{m}(x)\vert w(t,x)\vert^{2}dx\] 
\[=-(1+t)(w_{t}(t,\cdot),w(t,\cdot)) + (u_{1},u_{0}) + \frac{1}{2}\Vert w(t,\cdot)\Vert^{2} + \frac{1}{2}\int_{0}^{t}\int_{{\bf R}^{n}}a_{m}(x)\vert w(s,x)\vert^{2}dxds\]
\begin{equation}
 + \frac{1}{2}\int_{{\bf R}^{n}}a_{m}(x)\vert u_{0}(x)\vert^{2}dx + \int_{0}^{t}(1+s)\Vert w_{s}(s,\cdot)\Vert^{2}ds.
\end{equation}
Now, by using the Cauchy-Schwarz inequality and Lemma 2.2 we can first estimate
\[-(1+t)(w_{t}(t,\cdot),w(t,\cdot)) \leq \frac{(1+t)}{2\varepsilon}\Vert w_{t}(t,\cdot)\Vert^{2} + \frac{\varepsilon}{2}(1+t)\Vert w(t,\cdot)\Vert^{2}\]
\[\leq \frac{1+t}{\varepsilon}E_{w}(t) + \frac{\varepsilon}{2V_{0}}(1+t)\int_{{\bf R}^{n}}a_{m}(x)\vert w(t,x)\vert^{2}dx\]
\begin{equation}
\leq \frac{I_{1}^{2}}{\varepsilon} + \frac{\varepsilon}{2V_{0}}(1+t)\int_{{\bf R}^{n}}a_{m}(x)\vert w(t,x)\vert^{2}dx.
\end{equation}
(2.16) and (2.17) imply
\[\frac{1}{2}\Vert u_{0}\Vert^{2} + \int_{0}^{t}(1+s)\Vert\nabla w(s,\cdot)\Vert^{2}ds + \frac{(1+t)}{2}(1- \frac{\varepsilon}{V_{0}})\int_{{\bf R}^{n}}a_{m}(x)\vert w(s,x)\vert^{2}dx\] 
\[\leq (u_{1},u_{0}) + \frac{1}{2}\int_{{\bf R}^{n}}a(x)\vert u_{0}(x)\vert^{2}dx + \frac{1}{2}\Vert w(t,\cdot)\Vert^{2} + \frac{1}{2}\int_{0}^{t}\int_{{\bf R}^{n}}a_{m}(x)\vert w(s,x)\vert^{2}dxds
\]
\begin{equation}
+ \frac{I_{1}^{2}}{\varepsilon} + \int_{0}^{t}(1+s)\Vert w_{s}(s,\cdot)\Vert^{2}ds.
 \end{equation}
On the other hand, it follows from (2.5) and Lemma 2.2 that
\[\int_{0}^{t}(1+s)\Vert w_{s}(s,\cdot)\Vert^{2}ds \leq \frac{1}{V_{0}}\int_{0}^{t}(1+s)\int_{{\bf R}^{n}}a_{m}(x)\vert w_{s}(s,x)\vert^{2}dxds\]
\[= -\frac{1}{V_{0}}\int_{0}^{t}(1+s)E'_{w}(s)ds = -\frac{1}{V_{0}}(1+t)E_{w}(t) + \frac{1}{V_{0}}E(0) + \frac{1}{V_{0}}\int_{0}^{t}E_{w}(s)ds\]
\begin{equation}
\leq \frac{1}{V_{0}}(E(0) + I_{1}^{2}).
\end{equation}
Let us finalize the proof of Lemma 2.3. To do so, we rely on the inequality:
\[\frac{d}{dt}\{(1+t)^{2}E_{w}(t)\} \leq 2(1+t)E_{w}(t),\]
so that 
\begin{equation}
(1+t)^{2}E_{w}(t) \leq E(0) + \int_{0}^{t}(1+s)\Vert w_{s}(s,\cdot)\Vert^{2}ds + \int_{0}^{t}(1+s)\Vert \nabla w(s,\cdot)\Vert^{2}ds
\end{equation}
Because of Lemma 2.1, (2.17) with small $\varepsilon > 0$ and (2.18) one has
\[(1+t)^{2}E_{w}(t) \leq E(0) + \frac{2}{V_{0}}(E(0) + I_{1}^{2}) + (u_{1},u_{0}) + \frac{1}{2}\int_{{\bf R}^{n}}a(x)\vert u_{0}(x)\vert^{2} dx+ \frac{C}{2}I_{0}^{2} +\frac{I_{1}^{2}}{\varepsilon} =: I_{2}^{2} .\]

Concerning the fast $L^{2}$-decay estimate, we use (2.5) and (2.17) with small $\varepsilon > 0$ to have
\[\frac{(1+t)}{2}(1- \frac{\varepsilon}{V_{0}})V_{0}\int_{{\bf R}^{n}}\vert w(t,x)\vert^{2}dx \leq \frac{(1+t)}{2}(1- \frac{\varepsilon}{V_{0}})\int_{{\bf R}^{n}}a_{m}(x)\vert w(t,x)\vert^{2}dx\] 
\[\leq (u_{1},u_{0}) + \frac{1}{2}\int_{{\bf R}^{n}}a(x)\vert u_{0}(x)\vert^{2}dx + \frac{1}{2}\Vert w(t,\cdot)\Vert^{2} + \frac{1}{2}\int_{0}^{t}\int_{{\bf R}^{n}}a_{m}(x)\vert w(s,x)\vert^{2}dxds
\]
\begin{equation}
 +\frac{I_{1}^{2}}{\varepsilon}+ \int_{0}^{t}(1+s)\Vert w_{s}(s,\cdot)\Vert^{2}ds.
\end{equation}
The result follows from (2.20) and (2.18) and Lemma 2.1 by choosing $\varepsilon > 0$ small enough.
$\hfill\Box$
\par
\vspace{0.2cm}

{\it Proof of Theorem 1.1.} From Lemma 2.3 we first notice that $\{u^{(m)}\}$ is a bounded sequence in $L^{\infty}(0,\infty;H^{1}({\bf R}^{n}))$ and in $L^{\infty}(0,\infty;L^{2}({\bf R}^{n}))$. Furthermore, $\{u^{(m)}_{t}\}$ is also a bounded sequence in $L^{\infty}(0,\infty;L^{2}({\bf R}^{n}))$. Therefore, there exist a subsequence $\{u^{(\mu)}\}$ of the original one $\{u^{(m)}\}$, and a function $u = u(t,x) \in L^{\infty}(0,\infty;H^{1}({\bf R}^{n}))$ satisfying $u_{t} \in L^{\infty}(0,\infty;L^{2}({\bf R}^{n}))$ such that
\begin{equation}
u^{(\mu)} \to u\quad (\textstyle{weakly*}) \quad \textstyle{in} \quad L^{\infty}(0,\infty;H^{1}({\bf R}^{n})) \quad (\mu \to \infty),
\end{equation}
\begin{equation}
u_{t}^{(\mu)} \to u_{t} \quad (\textstyle{weakly*}) \quad \textstyle{in} \quad L^{\infty}(0,\infty;L^{2}({\bf R}^{n})) \quad (\mu \to \infty),
\end{equation}
\begin{equation}
u^{(\mu)} \to u \quad (\textstyle{weakly*}) \quad \textstyle{in} \quad L^{\infty}(0,\infty;L^{2}({\bf R}^{n})) \quad (\mu \to \infty).
\end{equation}
 
By multiplying both sides of the approximated equation (2.1) with $m$ replaced by $\mu$ the test function $\phi \in C_{0}^{\infty}([0,\infty)\times{\bf R}^{n})$, one can get the following weak form of the problem (2.1)-(2.2) with the help of the integration by parts:
\[\int_{0}^{\infty}\int_{{\bf R}^{n}}u^{(\mu)}(t,x)(\phi_{tt}(t,x)-\Delta\phi(t,x) - a_{\mu}(x)\phi_{t}(t,x))dx dt\] 
\begin{equation}
= \int_{{\bf R}^{n}}u_{1}(x)\phi(0,x)dx - \int_{{\bf R}^{n}}u_{0}(x)\phi_{t}(0,x)dx+ \int_{{\bf R}^{n}}a_{\mu}(x)u_{0}(x)\phi(0,x)dx.
\end{equation}

Now, it follows from (2.23) that as $\mu \to \infty$,
\begin{equation}
\int_{0}^{\infty}\int_{{\bf R}^{n}}u^{(\mu)}(t,x)(\phi_{tt}(t,x)-\Delta\phi(t,x))dx dt \to \int_{0}^{\infty}\int_{{\bf R}^{n}}u(t,x)(\phi_{tt}(t,x)-\Delta\phi(t,x))dx dt.
\end{equation}
Furthermore, because of the Lebesgue dominated convergence theorem one can get
\begin{equation}
\int_{{\bf R}^{n}}a_{\mu}(x)u_{0}(x)\phi(0,x)dx \to \int_{{\bf R}^{n}}a(x)u_{0}(x)\phi(0,x)dx \quad  (\mu \to \infty).
\end{equation}
On the other hand, for each fixed $\phi \in C_{0}^{\infty}([0,\infty)\times{\bf R}^{n})$, if we take $\mu \in {\bf N}$ large enough, then it follows from the compact support condition on $\phi(t,x)$ that
\[\int_{0}^{\infty}\int_{{\bf R}^{n}}u^{(\mu)}(t,x)a_{\mu}(x)\phi_{t}(t,x)dx dt = \int_{0}^{\infty}\int_{{\bf R}^{n}}u^{(\mu)}(t,x)a(x)\phi_{t}(t,x)dx dt.\] 
So, since $a(\cdot)\phi_{t}(t,\cdot) \in L^{1}([0,\infty);L^{2}({\bf R}^{n}))$, it follows from (2.23) that
\begin{equation}
\int_{0}^{\infty}\int_{{\bf R}^{n}}u^{(\mu)}(t,x)a_{\mu}(x)\phi_{t}(t,x)dx dt \to \int_{0}^{\infty}\int_{{\bf R}^{n}}u(t,x)a(x)\phi_{t}(t,x)dx dt \quad (\mu \to \infty).
\end{equation} 
Therefore, by taking $\mu \to \infty$ in (2.24), by means of (2.25)-(2.27) one can check that the limit function $u(t,x)$ is the weak solution to problem (1.1)-(1.2).\\

Finally, let us check decay estimates of the energy and $L^{2}$-norm of solutions such that
\begin{equation}
E_{u}(t) \leq C(1+t)^{-2}, \quad \Vert u(t,\cdot)\Vert^{2} \leq C(1+t)^{-1}
\end{equation}
with some constant $C > 0$. \\

For this end, one first remarks that for each $\phi \in L^{1}(0,\infty; C_{0}^{\infty}({\bf R}^{n}))$, it follows from the integration by parts that
\[\int_{0}^{\infty}(\frac{\partial u^{(\mu)}}{\partial x_{j}}(t,\cdot),\phi(t,\cdot))dt = -\int_{0}^{\infty}(u^{(\mu)}(t,\cdot), \frac{\partial \phi}{\partial x_{j}}(t,\cdot))dt,\]
so that we can have 
\[\lim_{\mu \to \infty}\int_{0}^{\infty}(\frac{\partial u^{(\mu)}}{\partial x_{j}}(t,\cdot),\phi(t,\cdot))dt = -\int_{0}^{\infty}(u(t,\cdot), \frac{\partial \phi}{\partial x_{j}}(t,\cdot))dt.\]
Since $u \in L^{\infty}([0,\infty);H^{1}({\bf R}^{n}))$, it follows from the integration by parts again one can get
\[\lim_{\mu \to \infty}\int_{0}^{\infty}(\frac{\partial u^{(\mu)}}{\partial x_{j}}(t,\cdot),\phi(t,\cdot))dt = \int_{0}^{\infty}(\frac{\partial u}{\partial x_{j}}(t,\cdot),\phi(t,\cdot))dt.\]
By density of $L^{1}(0,\infty; C_{0}^{\infty}({\bf R}^{n}))$ into $L^{1}(0,\infty;L^{2}({\bf R}^{n}))$, for each $j = 1,2,\cdots, n$ it is true that
\begin{equation}
\frac{\partial u^{(\mu)}}{\partial x_{j}} \to \frac{\partial u}{\partial x_{j}} \quad (\textstyle{weakly*}) \quad \textstyle{in} \quad L^{\infty}(0,\infty;L^{2}({\bf R}^{n})) \quad (\mu \to \infty).
\end{equation}
Now, let us finalize the proof of Theorem 1.1. First of all, we prepare the following basic lemma.
\begin{lem}\,Assume that a sequence $\{v_{m}\} \subset L^{\infty}(0,\infty;L^{2}({\bf R}^{n})$ satisfies 
\[v_{m} \to v \quad (\textstyle{weakly*}) \quad \textstyle{in} \quad L^{\infty}(0,\infty;L^{2}({\bf R}^{n})) \quad (m \to \infty),\]
for some $v \in L^{\infty}(0,\infty;L^{2}({\bf R}^{n}))$, and
\[\Vert v_{m}(t,\cdot)\Vert \leq C(1+t)^{-\gamma}\]
with some constants $C > 0$, and $\gamma > 0$. Then, it is also true that
\[\Vert v(t,\cdot)\Vert \leq C(1+t)^{-\gamma}.\]
\end{lem} 
{\it Proof.}\,Take $\psi \in C_{0}([0,\infty);L^{2}({\bf R}^{n}))$, and set $w_{m}(t,x) := (1+t)^{\gamma}v_{m}(t,x)$. Then, since\\
$(1+t)^{\gamma}\psi(t,x) \in L^{1}(0,\infty;L^{2}({\bf R}^{n}))$, by assumption it follows that
\[\lim_{m \to \infty}\int_{0}^{\infty}\int_{{\bf R}^{n}}w_{m}(t,x)\psi(t,x)dxdt = \lim_{m \to \infty}\int_{0}^{\infty}\int_{{\bf R}^{n}}v_{m}(t,x)\left((1+t)^{\gamma}\psi(t,x)\right)dxdt\]
\[= \int_{0}^{\infty}\int_{{\bf R}^{n}}v(t,x)\left((1+t)^{\gamma}\psi(t,x)\right)dxdt = \int_{0}^{\infty}\int_{{\bf R}^{n}}\left((1+t)^{\gamma}v(t,x)\right)\psi(t,x)dxdt.\]
Because of the density of $C_{0}([0,\infty);L^{2}({\bf R}^{n}))$ into $L^{1}(0,\infty;L^{2}({\bf R}^{n}))$ again (cf. Miyadera \cite[Theorem 15.3]{Miya}), we have
\[w_{m} = (1+t)^{\gamma}v_{m} \to (1+t)^{\gamma}v\quad (\textstyle{weakly*}) \quad \textstyle{in} \quad L^{\infty}(0,\infty;L^{2}({\bf R}^{n})) \quad (m \to \infty),\]
so that it follows that
\[\Vert (1+t)^{\gamma}v(t,\cdot)\Vert \leq \liminf_{m \to \infty}\Vert (1+t)^{\gamma}v_{m}\Vert_{L^{\infty}(0,\infty;L^{2}({\bf R}^{n}))} \leq C\]
for each $t \geq 0$. This implies the desired statement.
$\hfill\Box$

{\it Proof of Theorem 1.1 completed.}\, It follows from Lemma 2.3 that
\[
\Vert\frac{\partial u^{(\mu)}(t,\cdot)}{\partial x_{j}}\Vert \leq \sqrt{2}I_{2}(1+t)^{-1}, \quad \Vert u^{(\mu)}_{t}(t,\cdot)\Vert \leq \sqrt{2}I_{2}(1+t)^{-1}.
\]
Thus, it follows from (2.29) and (2.22) and Lemma 2.4 that
\begin{equation}
\Vert\frac{\partial u(t,\cdot)}{\partial x_{j}}\Vert \leq \sqrt{2}I_{2}(1+t)^{-1},
\end{equation}
\begin{equation}
\Vert u_{t}(t,\cdot)\Vert \leq \sqrt{2}I_{2}(1+t)^{-1}.
\end{equation}
Furthermore, from Lemma 2.3 one also has
\[\Vert u^{(\mu)}(t,\cdot)\Vert \leq I_{3}(1+t)^{-1/2}\quad (t \geq 0).\]
Therefore, (2.23) and Lemma 2.4 imply 
\begin{equation}
\Vert u(t,\cdot)\Vert \leq I_{3}(1+t)^{-1/2},\quad (t \geq 0).
\end{equation}
(2.30)-(2.32) imply the desired estimates (2.28) of Theorem 1.1. Note that all quantities $I_{j}$ ($j = 0,1,2,3$) can be absorbed into $CI_{00}$ defined in Theorem 1.1 with some constant $C > 0$. In this connection, we should remark the following relation because of the Cauchy-Schwarz inequality:
\[\int_{{\bf R}^{n}}a(x)\vert u_{0}(x)\vert^{2}dx \leq \sqrt{\int_{{\bf R}^{n}}a(x)^{2}\vert u_{0}(x)\vert^{2}dx}\sqrt{\int_{{\bf R}^{n}}\vert u_{0}(x)\vert^{2}dx}.\]
A uniqueness argument is standard, so we shall omit its detail.
$\hfill\Box$

\section{Remark on the low dimensional case.}

Let us give some remarks on the low dimensional case and several open problems.\\

{\rm (I)}\,When one gets the result for $n = 2$, we may use (2) of Proposition 2.1 with $\theta = 1$ and $\gamma \in (0,1]$ to get the similar estimate to (2.9), which is an essential part of proof. But, in this case we have to assume a stronger assumption such that
\begin{equation}
\int_{{\bf R}^{2}}(a(x)u_{0}(x) + u_{1}(x))dx = 0.
\end{equation}
Unfortunately, this assumption (3.1) is not hereditary to the approximate solution, i.e., 
$$\displaystyle{\int_{{\bf R}^{2}}}(a_{m}(x)u_{0}(x) + u_{1}(x))dx = 0$$
is not necessarily true. So, it is completely open to get the similar result in the low dimensional case (i.e., $n = 1,2$). \\

{\rm (II)}\,One can not treat the associated nonlinear equation such that
\begin{equation}
u_{tt}(t,x) -\Delta u(t,x) + a(x)u_{t}(t,x) = f(u(t,x)).\ \ \
\end{equation}
This is because one encounters lack of some compactness argument when we want to get the limit such that (see Lions \cite[(1.42)]{L})
\[f(u^{(m)}) \to f(u)\quad \textstyle{weakly*}\quad \textstyle{in}\quad L^{\infty}([0,\infty);X)\]
 (as $m \to \infty$), where $X$ is a Banach space.\\
In the unbounded coefficient case there are still many obstacles to be overcome. In this sense, the damped wave equation with unbounded coefficient and non-compactly supported initial data seem to be quite difficult to be treated.\\

{\rm (III)}\,We have another idea to deal with the two dimensional case, that is, the idea to rely on the well-known inequality:
\begin{equation}
\int_{{\bf R}^{2}}\frac{\vert\hat{f}(\xi)\vert^{2}}{\vert\xi\vert^{2}}d\xi \leq C\Vert f\Vert_{{\cal H}^{1}}^{2},
\end{equation}
in place of (2) of Proposition 2.1, where ${\cal H}^{1}({\bf R}^{2})$ is the so called Hardy space. Unfortunately, in this case we also encounter the same problem as the heredity of the property
$a(\cdot)u_{0} \in {\cal H}^{1}({\bf R}^{2})$ to $a_{m}(\cdot)u_{0} \in {\cal H}^{1}({\bf R}^{2})$.

\section{Appendix.}

In this section, let us check the density of $L^{1}(0,\infty; C_{0}^{\infty}({\bf R}^{n}))$ into $L^{1}(0,\infty;L^{2}({\bf R}^{n}))$.\\

Take $f \in L^{1}(0,\infty;L^{2}({\bf R}^{n}))$ and for each $L > 0$ set
\[ f_{L}(t,x) = \left\{
    \begin{array}{ll}
      \displaystyle{f(t,x)}&
           \qquad (\vert x\vert \leq L) \\[0.2cm]
      \displaystyle{0}& \qquad (\vert x\vert > L).
     \end{array} \right. \]
Next set
    
\[ h(t) = \left\{
    \begin{array}{ll}
      \displaystyle{e^{-\frac{1}{1-t^{2}}}}&
           \qquad (\vert t\vert < 1) \\[0.2cm]
      \displaystyle{0}& \qquad (\vert t\vert \geq 1),
     \end{array} \right. \]
and $\rho(x) := h(\vert x\vert^{2})\left(\displaystyle{\int_{{\bf R}^{n}}}h(\vert x\vert^{2})dx\right)^{-1}$, and for $\varepsilon > 0$ define
\[\rho_{\varepsilon}(x) := \frac{1}{\varepsilon^{n}}\rho(\frac{x}{\varepsilon}).\] 
Under these preparations, we define an approximating sequence of the original function $f$ by
\[f_{\varepsilon,L}(t,x) := (f_{L}(t,\cdot)*\rho_{\varepsilon})(x).\]
Then it is standard to check that $f_{\varepsilon,L}(t,\cdot) \in C_{0}^{\infty}({\bf R}^{n})$ for each $t \geq 0$ and $L > 0$, and
\[\Vert f_{\varepsilon,L}(t,\cdot)\Vert \leq \Vert f_{L}(t,\cdot)\Vert \leq \Vert f(t,\cdot)\Vert,\]
so that $f_{\varepsilon,L} \in L^{1}(0,\infty;L^{2}({\bf R}^{n}))$ for each $L > 0$ and $\varepsilon > 0$. Furthermore,  for each $L > 0$ and $t \geq 0$ it is known that 
\[\Vert f_{\varepsilon,L}(t,\cdot) - f_{L}(t,\cdot)\Vert \to 0\quad (\varepsilon \to 0).\]

Now, for an arbitrary fixed $\eta > 0$, choose $M > 0$ so large such that
\begin{equation}
\Vert f_{M} - f\Vert_{X} < \frac{\eta}{2},
\end{equation}
where $X := L^{1}(0,\infty;L^{2}({\bf R}^{n}))$, and $\Vert\cdot\Vert_{X}$ is the standard norm of $X$:
\[\Vert g\Vert_{X} := \int_{0}^{\infty}\Vert g(t,\cdot)\Vert dt.\]
Next, for such a fixed $M > 0$, by taking $\varepsilon > 0$ sufficiently small, if one applies the Lebesgue convergence theorem, one can get
\begin{equation}
\Vert f_{\varepsilon,M}-f_{M}\Vert_{X} < \frac{\eta}{2}.
\end{equation}
Therefore, it follows from (4.1) and (4.2) with such large $M > 0$ and small $\varepsilon > 0$ that
\[\Vert f_{\varepsilon,M} - f\Vert_{X} \leq \Vert f_{\varepsilon,M} - f_{M}\Vert_{X} + \Vert f_{M} - f\Vert_{X} < \frac{\eta}{2} + \frac{\eta}{2} = \eta,\]
which implies the density of $L^{1}(0,\infty; C_{0}^{\infty}({\bf R}^{n}))$ into $L^{1}(0,\infty;L^{2}({\bf R}^{n}))$.
$\hfill\Box$

\par
\vspace{0.5cm}
\noindent{\em Acknowledgement.}
\smallskip
The work of the first author (R. IKEHATA) was supported in part by Grant-in-Aid for Scientific Research 15K04958  of JSPS. 
The work of the second author (H. TAKEDA) was supported in part by Grant-in-Aid for Young Scientists (B)15K17581 of JSPS.


\end{document}